\documentclass[12pt]{article}

\usepackage{sectsty}
\usepackage{fullpage}
\usepackage[export]{adjustbox}
\usepackage{graphicx}

\usepackage{cancel,enumerate}
\usepackage[affil-it]{authblk}

\usepackage{amsmath}
\usepackage{amstext, amssymb}
\usepackage{eufrak}
\usepackage[mathscr]{euscript}

\usepackage{tikz}[2013/12/13]
\usetikzlibrary{calc}
\usetikzlibrary{matrix}
\definecolor{dark}{RGB}{0, 133, 202 }
\definecolor{light}{RGB}{112, 155, 230}

\newenvironment{proof}[1][Proof]{\noindent\textbf{#1.} }{\ \rule{0.5em}{0.5em}}

\newcommand{\E}{{\rm E}}
\newcommand{\hP}{\widehat{P}}
\newcommand{\hMu}{\widehat{\mu}}
\newcommand{\oX}{\overline{X}}
\newcommand{\hTh}{\widehat{\theta}}

\DeclareMathOperator{\NB}{\rm NB}
\newcommand{\st}{\,\vert\,}
\DeclareMathOperator{\Var}{Var}
\DeclareMathOperator{\Cov}{Cov}

\begin{document}

\title{Confidence Regions for Parameters of Negative Binomial Distribution}
\author[1]{Emmanuel Nkingi and Jan Vrbik}
\affil[1]{Department of Mathematics and Statistics\\ Brock University, Canada}
\date{\today}
\maketitle

\begin{abstract}
We describe a general method for the construction of a confidence region for the two parameters of the Negative Binomial Distribution. This is achieved by expanding the sampling distribution of Method-of-Moments estimators, using the Central Limit Theorem.
\end{abstract}

\section{\protect\bigskip Introduction}

Applied to a wide range of fields such as Biology, Medical Sciences, Communications, and Insurance, the Negative Binomial Distribution (NBD) has proven to be one of the most useful discrete distributions. As a result of its frequent application, an increasing number of papers on parameter estimation have appeared, and continue to appear, in literature (Fisher, 1941 and 1953; Gurland, 1959; Johnson and Kotz, 1969; Clark and Perry 1989; Chow, 1990; Shishebor and Towhidi, 2003; Famoye, 2011, among others).

In this article, we contribute further to this development by showing how to build an accurate confidence region for the unknown parameters of NBD.  A frequently used method for parameter estimation has been the Method of Moments. As such, our confidence region will be constructed based on this approach. There have been several choices of parameterization; in the field of Biology some have focused their attention on the dispersion parameter $\alpha$, and $p$; while others have used the reciprocal of $\alpha$. We use the expected value $\mu $ (location parameter) and shape parameter $P=\frac1{p}-1$ motivated by Johnson (Johnson et al., p.131).

\section{Background}

The NBD can be defined in two ways: The first version of the negative binomial distribution, first formulated by Montmort in 1714, is an
extension of the geometric distribution, where $X$ counts the number of  \emph{trials} until (and including) the $k^{\rm th}$ success.  In this version,
\begin{align}
P(X_{\rm old}=x\st p,k) =\dbinom{x-1}{k-1}p^{k}q^{x-k}; && x=k,k+1,\ldots \label{1}
\end{align}
Here $k$ is a positive integer; $0<p<1;$ and $q=(1-p)$.

The second version, which we consider for the purpose of this paper, counts the number of \emph{failures} before achieving the $k^{\rm th}$ success.  In this version,
\begin{align}
P(X_{\rm new}=x\st p,\alpha ) = \dbinom{x+\alpha -1}{x}p^{\alpha}q^{x}; && x=0,1,2,\ldots \label{2}
\end{align}

This version allows us to extend the definition of the negative binomial distribution to a positive \emph{real} parameter $\alpha $ (which was called $k$ in the previous distribution)$.$  Although it is impossible to visualize a non-integer number of ``successes'', we can still formally define the distribution through its probability mass function.

This form for dealing with NBD which has been used by our predecessors, will be inconvenient to manipulate in the manner desired for this paper; so, instead of using $\eqref2$ in its current form, we use the following parametrization
\begin{equation}
P(X=x \st \mu ,P) =\dbinom{x+\frac{\mu }{P}-1}{x}\left( \frac1{1+P}\right)^{\frac{\mu}{P}}\left(\frac{P}{1+P}\right)^{x};\qquad x=0,1,2,\ldots \label{new}
\end{equation}%
where as already mentioned in the introduction, 
\begin{align*}
p=\frac1{1+P} &&{\rm and}&& \alpha =\frac{\mu}{P}.
\end{align*}
The symbolic notation we use for this version will be $\mathrm{NB}( \mu,P)$.

The mean and the variance of $X$ are:
\begin{equation}
E(X) =\mu  \label{4}
\end{equation}
and
\begin{equation}
\sigma^2=\mu ( 1+P).  \label{5}
\end{equation}%
Thus, the variance is always larger than the mean for the NBD.  This property is at times referred to as \emph{over-dispersion}.

Note that the Poisson distribution is a limiting case of the $\NB( \mu ,P)$ when $P \to 0$ while keeping $\mu$ fixed.
The probability mass function of the Poisson distribution is,
\begin{equation}
\Pr[X=i]=\frac{\mu ^{i}}{i!}e^{-\mu }~\text{ where }~i=0,1,2,\ldots \label{8}
\end{equation}

\begin{proof}
The PGF (probability generating function) of the negative binomial distribution is given by:
\begin{equation}
\left( \frac1{1+P-Pz} \right)^{\frac{\mu }{P}}\   \label{91}
\end{equation}
as $P\to 0$, then the limit of \eqref{91} is
\begin{equation}
\frac1{e^{\mu (1-z)}}=e^{\mu (z-1)}  \label{12}
\end{equation}
which we recognize as the Poisson PGF.
\end{proof}

\section{Choosing Estimators}

An efficient method for estimating the parameters of most distributions is a procedure known as the Maximum Likelihood. \ It is well known that Maximum Likelihood Estimators (MLEs) are preferred over Method of Moment Estimators (MMEs) in the case of the NBD. However, when it comes to constructing confidence regions, MLEs, are rather difficult to deal with because they lead to a complicated function \eqref{70} for which we do not have an analytic expected value. 
\begin{equation}
\Var( \Psi ( \alpha +X) ) =\\E(\Psi^{^{\prime }}( \alpha ) -\Psi ^{^{\prime }}( \alpha+X) )  \label{70}
\end{equation}
where $\Psi (\alpha )$ is the digamma fun\bigskip ction,\ defined by:%
$$
\Psi (z)=\frac{\partial \log \Gamma (z)}{\partial z}=\frac{\Gamma ^{^{\prime}}(z)}{\Gamma (z)}~\text{ where }~z>0. 
$$
Using MLEs then requires performing a non-trivial numerical procedure for computation of this expected value.  For this reason, we have decided to go with MMEs instead, and the choice is justified by our results.  Moreover, from a practical point of view, our method is simple and can be
applied easily by anyone in the field of Statistics.

Investigating MLEs however, gave us a guideline as to how to deal with MMEs when faced with the condition of $(s^2<\oX)$.
Where
$$
s^2=\frac1{n}\sum_{i=1}^{n}( x_{i}-\oX) \text{.} 
$$
Unlike the MME which gives us a non-sensical answer (i.e. $P<0$), the MLE technique tells us that we have reached the Poisson limit. Being guided by the MLEs then, we take $P<0$ to be an indication that we have to switch to Poisson distribution, and set $P=0$.

\section{Method of Moments}

One of the most straightforward methods for estimating parameters of the Negative Binomial and other distributions is equating the first few sample moments to the corresponding theoretical moments, and solving the resulting equations for the unknown parameters.  This is known as the method of moments (MME).  The number of moments used is the same as the number of unknown parameters.

In the case of the Negative Binomial Distribution $($\ref2$)$ the equations for the first two moments are 
\begin{align} 
\hMu=\oX && \hMu( 1+\hP)=s^2  \label{13}
\end{align}
where $\oX$, and $s^2$ are the sample mean and sample variance, respectively.

The MME estimators of $\mu$ and $P$ are then equal to:
\begin{equation}
\widehat{\mu}=\oX  \label{14}
\end{equation}
and
\begin{equation}
\hP=\frac{s^2-\oX}{\oX} . \label{15}
\end{equation}

From \eqref{15} we see that we end up with a negative estimate $\hP$ whenever we are faced with a sample which has $s^2<\oX$. This again implies the distribution is Poisson.  To investigate the probability of this happening, we have utilized a Monte Carlo simulation.
\begin{table}
\begin{eqnarray*}
&&
\begin{tabular}{|l|l|l|l|l|}
\hline
\multicolumn{5}{|l|}{$n=30$} \\ \hline
$\mu \setminus P$ & $0.1$ & $0.3$ & $1$ & $10$ \\ \hline
$0.1$ & $79\%$ & $69\%$ & $54\%$ & $65\%$ \\ \hline
$0.3$ & $54\%$ & $36\%$ & $16\%$ & $25\%$ \\ \hline
$1$ & $49\%$ & $28\%$ & $0.05\%$ & $0.01\%$ \\ \hline
$10$ & $45\%$ & $23\%$ & $0.01\%$ & $0.00\%$ \\ \hline
\end{tabular}
\\
&&
\begin{tabular}{|l|l|l|l|l|}
\hline
\multicolumn{5}{|l|}{$n=50$} \\ \hline
$\mu \setminus P$ & $0.1$ & $0.3$ & $1$ & $10$ \\ \hline
$0.1$ & $67\%$ & $52\%$ & $36\%$ & $47\%$ \\ \hline
$0.3$ & $45\%$ & $24\%$ & $0.06\%$ & $0.1\%$ \\ \hline
$1$ & $42\%$ & $18\%$ & $0.01\%$ & $0.00\%$ \\ \hline
$10$ & $40\%$ & $14\%$ & $0.00\%$ & $0.00\%$ \\ \hline
\end{tabular}
\\
&&
\begin{tabular}{|l|l|l|l|l|}
\hline
\multicolumn{5}{|l|}{$n=100$} \\ \hline
$\mu \setminus P$ & $0.1$ & $0.3$ & $1$ & $10$ \\ \hline
$0.1$ & $46\%$ & $26\%$ & $0.12\%$ & $22\%$ \\ \hline
$0.3$ & $35\%$ & $11\%$ & $0.01\%$ & $0.01\%$ \\ \hline
$1$ & $32\%$ & $0.06\%$ & $0.0\%$ & $0.0\%$ \\ \hline
$10$ & $30$ & $0.04\%$ & $0.0\%$ & $0.0\%$ \\ \hline
\end{tabular}
\end{eqnarray*}
\caption{Occurrences of $s^2<\oX$ out of 10,000 random sample.}
\end{table}

We are quoting only the digits that are reliable.  A look at the above tables reveals, as expected, a high percentage of occurrences of $s^2<\oX$ when both parameters are small, while the percentage of $s^2<\oX$ begins to decrease as $\mu $ and $P$ get larger. \ If $P$ $\,$is much smaller than $\mu $, then $\sigma^2=\mu +P\approx \mu $; so, it is likely in this case to obtain samples with $s^2<\oX.$

Now that we have our estimators, we want to present them in the form of a confidence region, which gives us a better idea of what the actual parameters can be.  Since we do not have the exact joint distribution of $\mu$ and $P$, we have to use the Central Limit Theorem.

\section{Sampling Distribution}

\begin{figure}[htbp]
\centering
\includegraphics[trim={0em 18em 49em -1em},clip,frame]{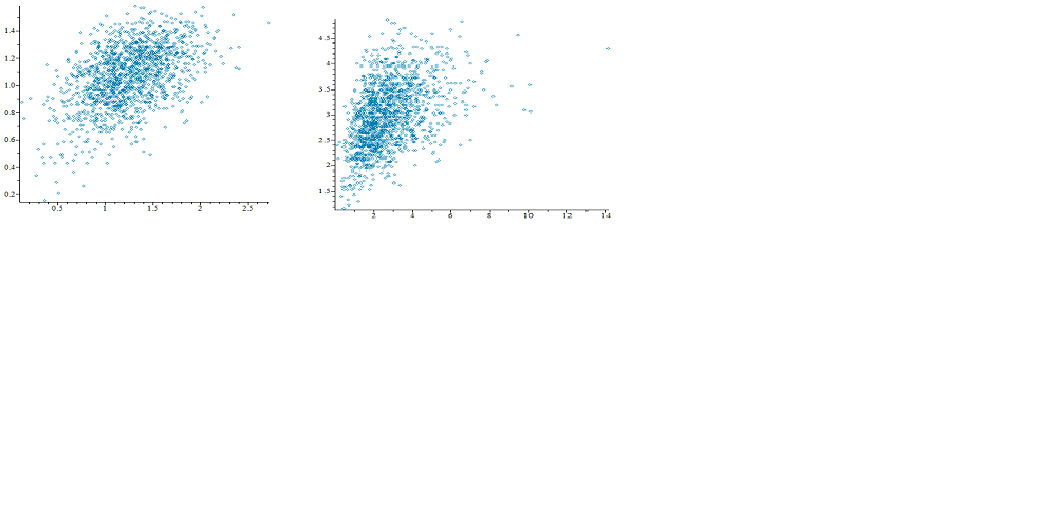}
\caption{A scatter plot of $\hMu$ and $\hP$ with $\mu =3$, $P=3$, and $n=30$.}
\label{534A}
\end{figure}

\begin{figure}[htbp]
\centering
\includegraphics[trim={20em 18em 27em -1em},clip,frame]{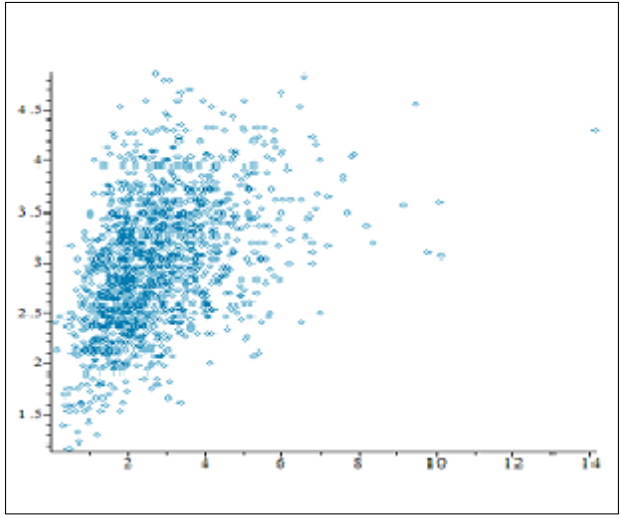}
\caption{A scatter plot of $\ln ( \hMu) $ and $\ln ( \hP+1) $ with $\mu =3$, $P=3$ and $n=30$.}
\label{534B}
\end{figure}

The Central Limit Theorem applies to all parameters based on sample \emph{means}.  Upon analyzing our $\hMu$ and $\hP$ distribution, we discovered that, for relatively small $n$ (i.e. $n=30$), their behaviour was far from a perfect bivariate normal distribution.  See Figure \ref{534A} and \ref{534B}.

We found that transforming our parameters to $\ln ( \mu ) $ and $\ln ( P+1) $ achieved a significant improvement towards normality and thus made our $CR$ more accurate.   See Figure \ref{534A}.

The corresponding method-of-moment estimators of ln$( \mu ) $ and $\ln( P+1) $ are
\begin{equation}
\ln ( \oX)  \label{20}
\end{equation}
and
\begin{equation}
\ln \left( \frac{\oX^2-\oX}{\oX} \right)  \label{21}
\end{equation}
respectively.

To construct CR for $\mu $ and $P$, we need to know the first few moments of the \emph{sampling distribution }of two estimators, \eqref{20} and \eqref{21}, namely the two means, two variances and covariance. We know that the expected value of any function of sample means can always be expanded, starting with the corresponding expected values.  Thus, to a sufficient approximation (ignoring the $\frac1{n}$ -- proportional correction) $\E(\ln ( \hMu) ) \approx \ln ( \mu ) $ and $E(\ln ( \hP+1) ) \approx \ln ( P+1)$.  What is left now is to find the formulas for the two variances and covariance.  To do that, we need to expand the two estimators at their respective means, up to the linear terms in $\oX$ and $\oX^2$, thus:
\begin{equation}
\ln ( \hMu) \approx \ln ( \mu ) +(\frac{\oX-\mu }{\mu }) + \cdots \label{22}
\end{equation}
and
\begin{equation}
\ln ( \hP+1) \approx \ln ( P+1) -\frac{1+P+2\mu }{\mu (1+P)}( \oX-\mu ) +\frac1{\mu (1+P)}(\oX^2-\E(X^2) )+ \cdots \label{23}
\end{equation}
where $\E(X^2) =\mu ( 1+\mu +P)$.

Replacing $z$ by $\exp ( t) $ in \eqref{91} yields the following MGF:
$$
( 1+P( 1-\exp( t) ) )^{-\frac{\mu }{P}} 
$$
which we expand to compute the first four simple moments of a single $X$:%
\begin{align*}
\E(X)   &= \mu  \\
\E(X^2) &= \mu ( 1+\mu +P)  \\
\E(X^3) &= \mu ( 1+3\mu +3P+3\mu P+\mu ^2+2P^2)  \\
\E(X^4) &= \mu ( 1+7\mu +7P+18\mu P+6\mu ^2+12P^2+6\mu^2P+11\mu P^2+\mu ^{3}+P^{3})
\end{align*}

Remembering that $\Var(\oX) $ is simply equal to $\frac{\Var( X) }{n},$ and similarly $\Var( \oX^2) = \frac{\Var( X^2) }{n}$, and $\Cov( \oX,\oX^2) =\frac{\Cov( X,X^2) }{n}$, we get:

\begin{align}
\Var( \oX)   &= \frac{\E(X^2) -\mu ^2}{n}= \frac{\mu ( 1+P) }{n}  \label{24} \\
\Var(\oX^2)  &= \frac{\E(X^{4}) -\E(X^2)^2}{n}=\frac{\mu(1+P) [ 2\mu ( 1+P) +P(2+3P) ] }{n}  \nonumber \\
\Cov( \oX,\, \oX^2) &= \frac{\E(X^{3})-\mu \E(X^2) }{n}=\frac{\mu ( 1+P) ( 1+2\mu+2P) }{n} . \nonumber
\end{align}

The variances and covariance of our estimators can now be easily computed
using \eqref{22}, \eqref{23}, and \eqref{24}:
\begin{align*}
\Var( \ln ( \hMu) ) &\equiv \frac{P+1}{\mu n}  \\
\Var( \ln ( \hP+1) ) &\equiv \frac{3P^2+2\mu P+2P+2\mu }{\mu ( 1+P) n}  \\
\Cov( \ln ( \hMu),\, \ln ( \hP+1)) &\equiv \frac{P}{\mu } .
\end{align*}

\section{Confidence Region}

Before discussing how to construct $CR$, we would like to make a more general statement about any two estimators which we are now going to call $\hTh_1$ and $\hTh_2$.  In the context of our discussion $\hTh_1=\ln ( \hMu) $ and $\hTh_2=\ln ( \hP+1)$.  Recall that any
general $\hTh$ can be standardized by:

\begin{equation}
Z_1=\frac{\hTh_1-\mu _1}{\frac{\sigma _1}{\sqrt{n}}}
\label{25}
\end{equation}%
and%
\begin{equation}
Z_2=\frac{\hTh_2-\mu_2}{\frac{\sigma_2}{\sqrt{n}}}
\label{26}
\end{equation}

The corresponding joint (bivariate) PDF is given by:
\begin{equation}
f(z_1,z_2)=\frac{\exp ( -\frac{z_1^2+z_2^2-2\rho z_1z_2}{2( 1-\rho ^2) }) }{2\pi \sqrt{1-\rho ^2}}  \label{27}
\end{equation}%
where $\mu_1$, $\sigma_1^2$, $\mu_{2}$, $\sigma_{2}^2$ are the means and variances of $\hTh_1$ and $\hTh_2$ respectively, and $\rho $ is their correlation coefficient.  Finding a contour of \eqref{27} still proves to be rather difficult, because of the correlation coefficient. \ To simplify matters we first transform $\hTh_1$ and $\hTh_2$ further to make them uncorrelated.

This can be achieved by
\begin{enumerate}
\item leaving the $\hTh_1$ unchanged, then 
\item we subtract a linear term $a\hTh_1$ from $\hTh_2$
\end{enumerate}
which yields: 
$$
\Cov( \hTh_2-a\hTh_1,\hTh_1)
=
\Cov( \hTh_1,\hTh_2)-a\Var( \hTh_1)
= 0 \implies a
=\frac{\Var( \hTh_1) }{\Cov( \hTh_1,\hTh_2) }  .
$$  
In our case $a=\frac{P}{1+P}$ implying that 
\begin{align}
\hTh_2-a\hTh_1 
&= \ln ( \hP+1)-\ln ( P+1) -\frac{P}{1+P}( \ln ( \hMu)-\ln ( \mu ) )  \label{50}
\Var( \hTh_2-a\hTh_1) 
&= ( \sigma_2^2-\rho ^2\sigma _2^2) .
\end{align}
This equals to $2( \frac{\mu +P}{n\mu }) $, which is the variance of \eqref{50}.

\subsection{New Standardization}

\begin{figure}[htbp]
\centering
\includegraphics[trim={0em 0em 0em 0em},clip,frame]{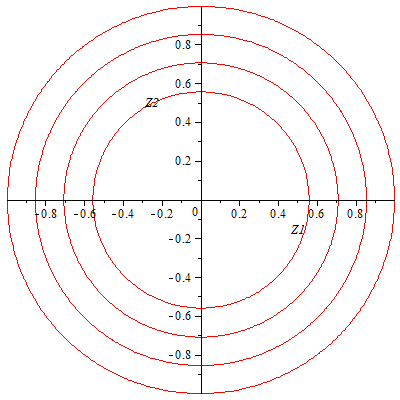}
\caption{Contour plot for \eqref{ThirtyNine}}
\label{540}
\end{figure}

We define $Z_1$ the same way as \eqref{25}
\begin{equation}
Z_1=\frac{\hTh_1-\mu _1}{\frac{\sigma _1}{\sqrt{n}}}
\label{28}
\end{equation}
but now $Z_2$ changes to
\begin{equation}
Z_2=\frac{\hTh_2-\mu _2-\frac{\rho \sigma _2}{\sigma_1}( \hTh_1-\mu _1) }{\sigma _2\sqrt{\frac{1-\rho ^2}{n}}} .  \label{29}
\end{equation}
This yields the desired PDF for building $CR$ namely
\begin{equation}
f(z_1,z_2)=\frac{\exp ( -\frac{z_1^2+z_2^2}2) }{2\pi }  \label{ThirtyNine}
\end{equation}
with the corresponding contour plot given in Figure \ref{540}.

Constructing a $100( 1-\delta )$ joint confidence region for any  $\theta_1$ and $\theta_2$ now amounts to solving 
\begin{equation}
\frac{\exp ( -\frac{z_1^2+z_2^2}2) }{2\pi }\geq c \label{32}
\end{equation}%
or simply 
\begin{equation}
z_1^2+z_2^2\preceq c_0,  \label{33}
\end{equation}
so that
\begin{equation}
\int \int_{z_1^2+z_2^2\preceq c_0}\frac{\exp ( -\frac{z_1^2+z_2^2}2) }{2\pi } ~ dz_1dz_2=1-\delta  .\label{34}
\end{equation}
This yields
\begin{equation}
c_0=-2\ln ( \delta ) . \label{35}
\end{equation}

\begin{proof}
Switching to polar coordinates, and solving for $c_0$ amounts to solving $%
r^2=c_0$, thus:
$$
\int_0^{2\pi }\int_0^{\sqrt{c_0}}\frac{\exp ( -\frac{r^2}2) }{2\pi } ~rdrd\vartheta 
= \int_0^{\frac{c_0}2}\exp (-u) ~du=1-\exp ( -\frac{c_0}2) 
= 1-\delta 
$$
implying \eqref{35}.
\end{proof}

In the context of our distribution, the $CR$ is given by:%
\begin{equation}
\frac{\left[ \ln ( \hMu) -\ln ( \mu ) \right] ^2}{\frac{1+P}{n\mu }}+\frac{\left[ \ln ( \hP+1) -\ln ( P+1) -\frac{P}{1+P}( \ln ( \hMu) -\ln ( \mu ) ) \right] ^2}{\frac{2( \mu
+P) }{n\mu }} \leq -2\ln ( \delta ) \label{40} 
\end{equation}
where $n$, $\delta$, $\widehat{\mu}$, and $\hP+1$ must be replaced by their numerical values: $1-\delta $ specifies the desired level of confidence.

\subsection{Example}

Consider a computer generated random independent sample of size $50$ from the NBD with $\mu =1$ and $P=1$ (assumed unknown and to be estimated), namely
\begin{align*}
0,1,2,0,0,0,0,3,0,0,1,0,1,1,2,1,0,0,0,3,0,0,0,3,0 \\
1,0,1,1,0,0,1,1,1,3,7,1,1,0,1,3,0,0,0,0,0,3,0,1,2
\end{align*}
Computing $\oX=\sum_{i=1}^{50}\frac{x_{i}}{n}=0.96$ and $\oX^2=\sum_{i=1}^{50}\frac{x_{i}^2}{n}=2.60$ allows us evaluate the MME expressions of \eqref{20} and \eqref{21}, getting
\begin{align}
\widehat{\mu}=0.960 &&\text{and}&& \hP+1=1.906 . \label{41}
\end{align}
To find (say) a $50\%,$ $80\%,$ and $95\%$ confidence regions for $\mu$ and $P$, one needs to solve \eqref{40} for each $\delta$.  Evaluating \eqref{40} with the estimates found in \eqref{41}, we get
\begin{equation}
\frac{\left[ \ln ( 0.960) -\ln ( \mu ) \right] ^2}{\frac{1+P}{50\mu }}+\frac{\left[ \ln ( 1.906) -\ln (
P+1) -\frac{P}{1+P}( \ln ( 0.960) -\ln ( \mu ) ) \right] ^2}{\frac{2( \mu +P) }{50\mu }} 
\leq -2\ln ( \delta ). \label{51} 
\end{equation}

\begin{figure}[htbp]
\centering
\includegraphics[trim={-1em 8em 42em -1em},clip,frame]{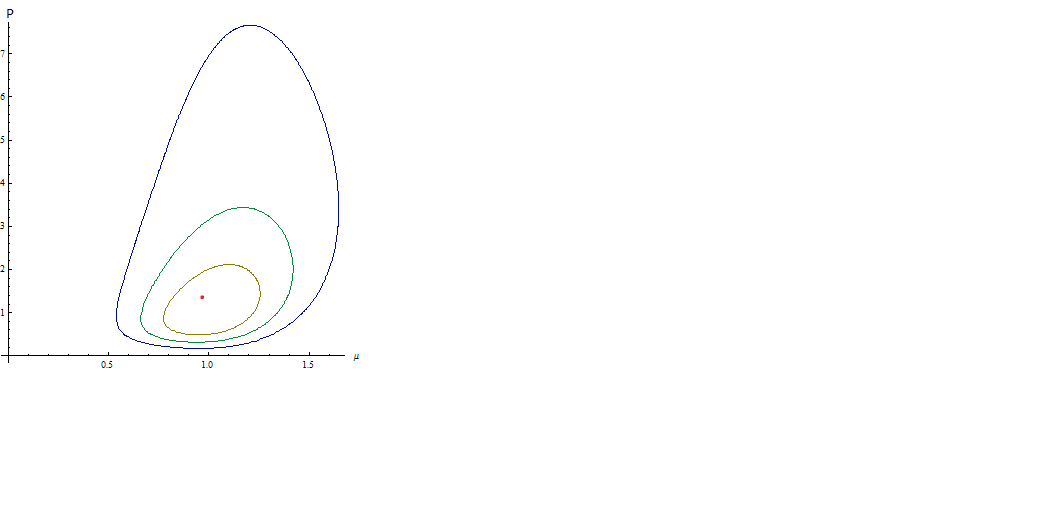}
\caption{Confidence regions for $\mu $ and $P$, with $\widehat{\mu}=0.960$ and $\hP=1.906$.}
\label{542}
\end{figure}

Even though \eqref{51} is a non-linear equation in terms of $\mu$ and $P$, plotting the actual contours is fairly simple with the help of a computer.  Thus, for a $95\%$, $80\%$, and $50\%$ confidence level, we Figure \ref{542}.

We find that all our confidence regions contain the ``true'' values used in our simulation.  We investigate the accuracy of this method in the following section.  Impressively, our method works even in the $s^2<\oX$ situation.

\subsection{Example}

Suppose now that our generated random sample of size 50 is drawn from a NBD with $\mu =3$ and $P=0.3$, and consists of:%
\begin{align*}
0,1,0,0,2,0,0,1,0,0,1,0,1,2,0,1,0,3,2,1,1,2,2,2,2 \\
2,1,0,0,0,3,1,0,1,0,0,1,2,0,0,0,1,3,1,1,1,0,1,0,0
\end{align*}

We once again compute $\oX=\sum_{i=1}^{50}\frac{x_{i}}{n}=2.98$ and $\oX^2=\sum_{i=1}^{50}\frac{x_{i}^2}{n}=2.880\overline{6}$ which allows us to find the MME using \eqref{20} and \eqref{21}
\begin{equation}
\hMu=2.98\text{ and }\hP=-0.03\overline{3}  \label{60}
\end{equation}

\begin{figure}[htbp]
\centering
\includegraphics[trim={-1em 10em 42em -1em},clip,frame]{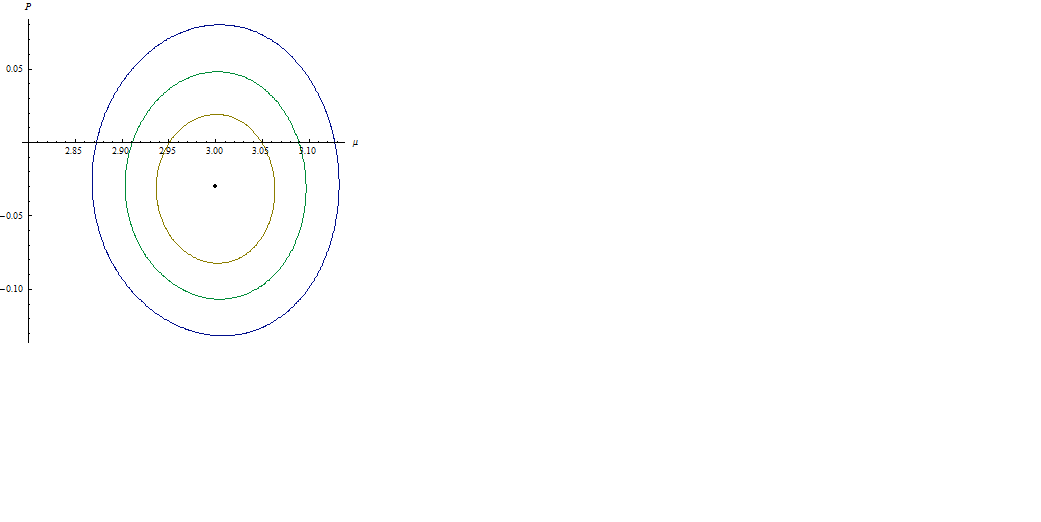}
\caption{Confidence regions for $\mu$ and $P$, with $\hMu=3$ and $\hP=-0.03\overline{3}$.}
\label{545}
\end{figure}

Clearly $P$ cannot be negative.  However, when the NBD distribution is close to the Poisson limit, the MME estimate may easily become negative, unlike the MLE estimate which becomes zero.  Even though we have a nonsensical estimate in $\hP$, we can still build an accurate confidence region using the same formulas.  Solving \eqref{40} with the same $\delta$ as the previous example, we get the contour plot of Figure \ref{545}.

\begin{figure}[htbp]
\centering
\includegraphics[trim={-1em 9em 42em -1em},clip,frame]{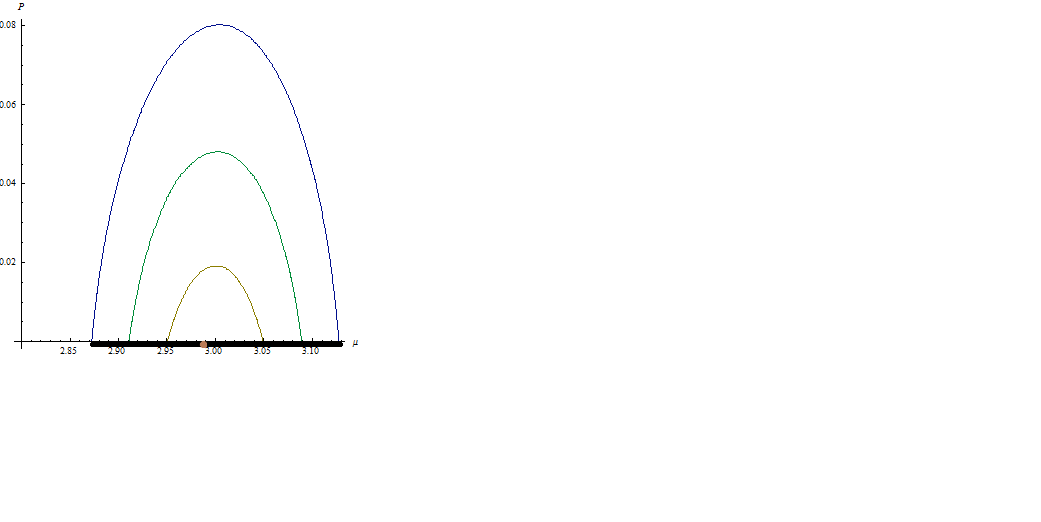}
\caption{Confidence regions for $\mu $ and $P$, with $\hMu=3$ and $\hP=-0.03\overline{3}$.}
\label{546}
\end{figure}

As seen in Figure \ref{545}, there is a strong indication that the distribution is Poisson. However, there is also the possibility that it is not.  The region under the $\mu$- axis is what we consider to imply Poisson. Conversely, the region above is Negative Binomial.  In the MLE case, the negative region is never observed; the entire non-physical region is clasped on the $\mu$-axis.  Using this `hint', we can present the same results in a new, more meaningful form of Figure \ref{546}.

\section{Monte Carlo Verification}

In order to compare how well our approximation matches the respective confidence levels, we generate a large collection of random independent samples (say $10,000$) from the NBD, using a specific choice of $\mu$ and $P$, and establish, for each of these samples, the percentage of the corresponding confidence regions which cover the two values $\mu$ and $P$.  When done correctly, the probability of this happening should be close to the chosen significant level. We then repeat this for as many different choices of $\mu$ and $P$ as feasible.

\begin{table}[htbp]
$$
\begin{tabular}{|l|l|lll|}
\hline
\multicolumn2{|l}{$100( 1-\delta ) \%$} & $50\%$ & $80\%$ & $95\%$ \\ \hline
$\mu \setminus P$ & $n$ & \multicolumn{3}{|c|}{$0.3$} \\ 
\multicolumn1{|c|}{} & $30$ & $53.35\%$ & $85.95\%$ & $95.55\%$ \\ 
$0.3$ & $50$ & $51.87\%$ & $83.44\%$ & $95.68\%$ \\ 
& $100$ & $52.94\%$ & $80.93\%$ & $95.38\%$ \\ \hline
& $30$ & $51.15\%$ & $81.05\%$ & $94.58\%$ \\ 
$1$ & $50$ & $50.24\%$ & $80.65\%$ & $95.71\%$ \\ 
& $100$ & $49.69\%$ & $80.67\%$ & $94.91\%$ \\ \hline
& $30$ & $48.58\%$ & $78.70\%$ & $94.25\%$ \\ 
$3$ & $50$ & $49.86\%$ & $79.71\%$ & $94.78\%$ \\ 
& $100$ & $50.15\%$ & $79.89\%$ & $95.29\%$ \\ \hline
\end{tabular}%
$$
$$
\begin{tabular}{|l|l|lll|}
\hline
\multicolumn2{|l}{$100( 1-\delta ) \%$} & $50\%$ & $80\%$ & $%
95\%$ \\ \hline
$\mu \setminus P$ & $n$ & \multicolumn{3}{|c|}{$0.3$} \\ 
\multicolumn1{|c|}{} & $30$ & $56.43\%$ & $86.20\%$ & $96.33\%$ \\ 
$0.3$ & $50$ & $56.55\%$ & $85.28\%$ & $95.59\%$ \\ 
& $100$ & $54.25\%$ & $83.27\%$ & $95.87\%$ \\ \hline
& $30$ & $52.27\%$ & $81.13\%$ & $95.67\%$ \\ 
$1$ & $50$ & $51.32\%$ & $81.33\%$ & $95.36\%$ \\ 
& $100$ & $50.90\%$ & $81.54\%$ & $95.34\%$ \\ \hline
& $30$ & $49.91\%$ & $79.86\%$ & $94.60\%$ \\ 
$3$ & $50$ & $50.05\%$ & $79.54\%$ & $94.51\%$ \\ 
& $100$ & $50.05\%$ & $79.60\%$ & $94.68\%$ \\ \hline
\end{tabular}%
$$
$$
\begin{tabular}{|l|l|lll|}
\hline
\multicolumn2{|l}{$100( 1-\delta ) \%$} & $50\%$ & $80\%$ & $%
95\%$ \\ \hline
$\mu \setminus P$ & $n$ & \multicolumn{3}{|c|}{$0.3$} \\ 
\multicolumn1{|c|}{} & $30$ & $59.33\%$ & $86.62\%$ & $96.58\%$ \\ 
$0.3$ & $50$ & $58.30\%$ & $86.23\%$ & $94.96\%$ \\ 
& $100$ & $57.06\%$ & $86.17\%$ & $96.78\%$ \\ \hline
& $30$ & $55.92\%$ & $85.31\%$ & $96.22\%$ \\ 
$1$ & $50$ & $55.20\%$ & $84.88\%$ & $96.64\%$ \\ 
& $100$ & $52.54\%$ & $82.99\%$ & $96.20\%$ \\ \hline
& $30$ & $53.13\%$ & $82.75\%$ & $95.38\%$ \\ 
$3$ & $50$ & $51.33\%$ & $81.75\%$ & $95.88\%$ \\ 
& $100$ & $50.65\%$ & $80.19\%$ & $95.22\%$ \\ \hline
\end{tabular}%
$$
\caption{Monte Carlo Verification}
\label{Table2}
\end{table}

As a result, we can observe a reasonably close match between the desired and
actually achieved value of $\delta $.

\section{Conclusion}

We have demonstrated how to find confidence regions for unknown parameters of the negative binomial distribution. Uniquely, the method described in this article works even in the case of over-dispersed samples (i.e. $s^2<\oX$).  In this case, we find that our confidence region split into two regions, one corresponding to the Poisson and the other to Negative Binomial distribution.

One could now consider the Edgeworth Expansion (improving over normal approximation), but, considering how good the results were in Table 2, and also the extra difficulties one runs into when lowering the sample size, it is something we saw as a direction of decreasing marginal returns.  We are, however, looking forward to contrasting our results with a future paper in which we will show how to build confidence regions for unknown parameters of the NBD using the MLE method.

\end{document}